\documentclass[a4paper,reqno]{amsart}

\usepackage[T1]{fontenc}
\usepackage[latin1]{inputenc}

\PII{}
\makeatletter
\def\@setcopyright{\@empty}
\makeatother

\newcommand{\seq}[3]{\{#1\}_{#2}^{#3}}

\newcommand{\numericset}[1]{\mathbb #1}
\newcommand{\numR}{\numericset R}

\theoremstyle{plain}
\newtheorem{theorem}{Theorem}[section]

\newtheorem{corollary}{Corollary}[section]

\theoremstyle{definition}

\theoremstyle{remark}

\newcounter{const}[section]
\numberwithin{const}{theorem}
\numberwithin{const}{lemma}
\numberwithin{const}{corollary}
\numberwithin{const}{example}

\makeatletter
\newcommand{\Cn}[1][]{%
  \stepcounter{const}C_{\theconst}%
  \@ifnotempty{#1}{\newcounter{#1}\setcounter{#1}{\arabic{const}}}}
\makeatother

\newcommand{\lastC}{C_{\theconst}}
\newcommand{\prevC}[1][1]{%
	{\countdef\n=255
	 \n=\theconst
	 \advance\n by-#1
	 C_{\number\n}}}

\numberwithin{equation}{section}

\renewcommand{\theconst}{\arabic{const}}

\begin{document}

\title[Some reverse $l_p$-type inequalities\dots]%
	{%
		Some reverse $l_p$-type inequalities
		involving certain quasi monotone sequences%
	}

\author[M.~K.\ Potapov]{Mikhail~K.\ Potapov}
\address{Mikhail~K.\ Potapov\\
	Department of Mechanics and Mathematics\\
	Moscow State University\\
	Moscow 117234\\
	Russia}
\email{mkpotapov@mail.ru}

\author[F.~M.\ Berisha]{Faton~M.\ Berisha}
\address{Faton~M.\ Berisha\\
	Faculty of Mathematics and Sciences\\
	University of Prishtina\\
	Nëna Terezë~5\\
	10000 Prishtina\\
	Kosovo%
}
\email{faton.berisha@uni-pr.edu}

\author[N.~Sh.\ Berisha]{Nimete~Sh.\ Berisha}
\address{Nimete~Sh.\ Berisha\\
	Faculty of Economics\\
	University of Prishtina%
}
\email{nimete.berisha@gmail.com}

\author[R.~Kadriu]{Reshad Kadriu}
\address{Reshad Kadriu\\
	College "Business"\\
	Prishtina%
}

\date{January 9, 2014}

\keywords{%
	$l_p$-type, Copson, Leindler, inequalities,
	quasi, lacunary, geometrically,
	monotone sequences%
}
\subjclass{Primary 47A30, Secondary 26D15.}
\date{}

\begin{abstract}
	In this paper,
	we give some $l_p$-type inequalities
	about sequences
	satisfying certain quasi monotone type properties.
	As special cases,
	reverse $l_p$-type inequalities
	for non-negative decreasing sequences are obtained.
	The inequalities are closely related to Copson's
	and Leindler's inequalities,
	but the sign of the inequalities is reversed.
\end{abstract}

\maketitle

\section{Introduction}

For non-negative number sequences the following,
classical inequalities of Hardy and Littlewood,
are well known~\cite[p.~255, Th~346]{hardy-l-p:inequalities}.

Let $\seq{b_\nu}{\nu=1}{\infty}$ be a sequence of non-negative numbers,
$\alpha>0$, $m$ and~$n$ positive integers such that $n<m$.
The following inequalities hold true:
\begin{gather}
	\sum_{\mu=n}^m \mu^{\alpha-1}
		\biggl(
			\sum_{\nu=\mu}^m b_\nu
		\biggr)^p
	\le C
	\sum_{\mu=n}^m \mu^{\alpha-1}
		(\mu b_\mu)^p,
	\label{eq:p-le-pge1}
	\\
	\sum_{\mu=n}^m \mu^{-\alpha-1}
		\biggl(
			\sum_{\nu=n}^\mu b_\nu
		\biggr)^p
	\le C
	\sum_{\mu=n}^m \mu^{-\alpha-1}
		(\mu b_\mu)^p
	\label{eq:m-le-pge1}
\end{gather}
for $p\ge1$;
and
\begin{gather}
	\sum_{\mu=n}^m \mu^{\alpha-1}
		\biggl(
			\sum_{\nu=\mu}^m b_\nu
		\biggr)^p
	\ge C
	\sum_{\mu=n}^m \mu^{\alpha-1}
		(\mu b_\mu)^p,
	\label{eq:p-ge-ple1}
	\\
	\sum_{\mu=n}^m \mu^{-\alpha-1}
		\biggl(
			\sum_{\nu=n}^\mu b_\nu
		\biggr)^p
	\ge C
	\sum_{\mu=n}^m \mu^{-\alpha-1}
		(\mu b_\mu)^p
	\label{eq:m-ge-ple1}
\end{gather}
for $0<p\le1$,
where positive constant~$C$
depends only on numbers~$\alpha$ and~$p$,
and does not depend on integers~$m$, $n$,
and the sequence $\seq{b_\nu}{\nu=1}\infty$.

Closely related to these inequalities
are classical Copson inequalities~\cite{copson:london-28},
Leindler's inequalities%
~\cite{leindler:acta-66, leindler:acta-70, leindler:jipam-00},
and those proved
or used in%
~\cite{potapov:math-72, potapov-b:publ-79,
	potapov-s-t:approx-13, gao:arxiv-1301%
}.

In the paper we prove some related inequalities
which involve non-negative sequences
satisfying certain
monotone-type properties.
As special cases,
inequalities converse to~\eqref{eq:p-le-pge1},
\eqref{eq:m-le-pge1}, \eqref{eq:p-ge-ple1}
and~\eqref{eq:m-ge-ple1}
for the case of non-negative monotone decreasing number sequences
are deduced.

In order to prove the inequalities
we need the following
\begin{theorem}\label{th:jensen}
	Let $\seq{b_\nu}{\nu=1}\infty$ be a sequence of non-negative numbers,
	$0<\alpha<\beta$,
	$m$ and~$n$ positive integers such that $n<m$.
	Then the following inequality holds
	$$
	\biggl(
		\sum_{\mu=n}^m b_\mu^\beta
	\biggr)^{1/\beta}
	\le
	\biggl(
		\sum_{\mu=n}^m b_\mu^\alpha
	\biggr)^{1/\alpha}.
	$$
\end{theorem}

The proof of the theorem
is due to Jensen~\cite[p.~28, Th.~19]{hardy-l-p:inequalities}.

We call a number sequence $\seq{a_\nu}{\nu=1}\infty$
non-negative monotone decreasing (or increasing),
and denote it by $a_\nu\downarrow$ (or $a_\nu\uparrow$),
if for each positive integer~$\nu$
the following conditions are satisfied
\begin{enumerate}
	\item $a_\nu\ge0$,
	\item $a_{\nu+1}\le a_\nu$
		\quad (or $a_{\nu+1}\ge a_\nu$, respectively).
\end{enumerate}

We call a number sequence $\seq{\lambda_\nu}{\nu=1}\infty$
quasi lacunary monotone
if it is a non-negative monotone sequence
(i.e., $\lambda_\nu\downarrow$ or $\lambda_\nu\uparrow$)
and there are positive constants~$K_1$ and~$K_2$
such that for each positive integer~$\nu$
the following condition is satisfied
\begin{displaymath}
	K_1\lambda_{2^\nu}
	\le\lambda_{2^{\nu+1}}
	\le K_2\lambda_{2^\nu}.
\end{displaymath}

Finally, we call a non-negative number sequence
$\seq{\lambda_\mu}{\mu=1}\infty$
quasi geometrically increasing
if there is a positive constant~$K$
such that for each positive integer~$m$
the following condition is satisfied
\begin{displaymath}
	\sum_{\mu=1}^m \lambda_\mu
	\le K\lambda_m.
\end{displaymath}

\section{Inequalities for quasi monotone sequences}

In this section we give several reverse inequalities of $l_p$-type
involving non-negative decreasing sequences,
quasi lacunary monotone sequences
or quasi geometrically increasing sequences. 

\begin{theorem}\label{th:lambda-ge-pge1}
	Let a sequence $\seq{a_\nu}{\nu=1}\infty$
	be such that $a_\nu\downarrow$,
	$\seq{\lambda_\mu}{\mu=1}\infty$
	and $\seq{\gamma_\nu}{\nu=1}\infty$
	be lacunary monotone sequences,
	$m$ and~$n$ positive integers such that $m\ge16n$.
	If $p\ge1$,
	then the following inequality holds
	\begin{equation}\label{eq:lambda-m-ge-pge1}
		\sum_{\mu=n}^m\lambda_\mu
			\biggl(
				\sum_{\nu=n}^\mu a_\nu\gamma_\nu
			\biggr)^p
		\ge C
		\sum_{\mu=4n}^m\lambda_\mu
			(\mu a_\mu\gamma_\mu)^p.
	\end{equation}
	If, in addition,
	$\seq{2^\mu\lambda_{2^\mu}}{\mu=1}\infty$
	is a quasi geometrically increasing sequence,
	then the following inequality holds
	\begin{equation}\label{eq:lambda-p-ge-pge1}
		\sum_{\mu=n}^m\lambda_\mu
			\biggl(
				\sum_{\nu=\mu}^m a_\nu\gamma_\nu
			\biggr)^p
		\ge C
		\sum_{\mu=8n}^m\lambda_\mu
			(\mu a_\mu\gamma_\mu)^p.
	\end{equation}
	Here and later on
	$C$ and~$C_i$ will denote positive constants
	depending only on~$p$
	and the sequences $\seq{\lambda_\mu}{\mu=1}\infty$
	and $\seq{\gamma_\nu}{\nu=1}\infty$,
	and not depending on~$m$, $n$
	and the sequence $\seq{a_\nu}{\nu=1}\infty$.
\end{theorem}
                           
\begin{theorem}\label{th:lambda-le-ple1}
	Let a sequence $\seq{a_\nu}{\nu=1}\infty$
	be such that $a_\nu\downarrow$,
	$\seq{\lambda_\mu}{\mu=1}\infty$
	and $\seq{\gamma_\nu}{\nu=1}\infty$
	be lacunary monotone sequences,
	$m$ and~$n$ positive integers such that $m\ge4n$.
	If $0<p\le1$,
	then the following inequality holds
	\begin{equation}\label{eq:lambda-m-le-ple1}
		\sum_{\mu=4n}^m\lambda_\mu
			\biggl(
				\sum_{\nu=4n}^\mu a_\nu\gamma_\nu
			\biggr)^p
		\le C
		\sum_{\mu=n}^m\lambda_\mu
			(\mu a_\mu\gamma_\mu)^p.
	\end{equation}
	If, in addition,
	$\seq{2^\mu\lambda_{2^\mu}}{\mu=1}\infty$
	is a quasi geometrically increasing sequence,
	then the following inequality holds
	\begin{equation}\label{eq:lambda-p-le-ple1}
		\sum_{\mu=4n}^m\lambda_\mu
			\biggl(
				\sum_{\nu=\mu}^m a_\nu\gamma_\nu
			\biggr)^p
		\le C
		\sum_{\mu=n}^m\lambda_\mu
			(\mu a_\mu\gamma_\mu)^p.
	\end{equation}
\end{theorem}

\section{Proof of Theorem~\ref{th:lambda-ge-pge1}}

We prove inequality~\eqref{eq:lambda-p-ge-pge1}.
For given~$n$ and~$m$
we choose positive integers~$N$ and~$M$
such that
$2^{N-1}<n\le2^N$
and $2^M\le m<2^{M+1}$.
Then the following inequality holds
\begin{displaymath}
	I=
	\sum_{\mu=n}^m\lambda_\mu
		\biggl(
			\sum_{\nu=\mu}^m a_\nu\gamma_\nu
		\biggr)^p
	\ge
	\sum_{\mu=2^N}^{2^M}\lambda_\mu
		\biggl(
			\sum_{\nu=\mu}^{2^M} a_\nu\gamma_\nu
		\biggr)^p.
\end{displaymath}
By splitting the first sum into blocks of length $2^i$,
we obtain
\begin{displaymath}
	I\ge
	\sum_{i=N+1}^M
		\sum_{\mu=2^{i-1}+1}^{2^i}
			\lambda_\mu
			\biggl(
				\sum_{\nu=\mu}^{2^M} a_\nu\gamma_\nu
			\biggr)^p.
\end{displaymath}
By bounding the third sum from below,
taking into account that $\seq{\lambda_\mu}{\mu=1}\infty$
is a quasi lacunary monotone sequence,
we have
\begin{displaymath}
	I\ge
	\sum_{i=N+1}^M
		\biggl(
			\sum_{\nu=2^i}^{2^M} a_\nu\gamma_\nu
		\biggr)^p
		\sum_{\mu=2^{i-1}+1}^{2^i}
			\lambda_\mu
	\ge
	\Cn
	\sum_{i=N+1}^M
		2^{i-1}\lambda_{2^{i-1}}
		\biggl(
			\sum_{\nu=2^i}^{2^M} a_\nu\gamma_\nu
		\biggr)^p,
\end{displaymath}
where $\lastC$
depends only on the sequence $\seq{\lambda_\mu}{\mu=1}\infty$.
Now, we split the second sum into blocks of length $2^{i-1}$,
remove the terms with index $i=M$,
and taking into consideration that $a_\nu\downarrow$
and $\seq{\lambda_\mu}{\mu=1}\infty$,
$\seq{\gamma_\nu}{\nu=1}\infty$
are lacunary monotone sequences,
we get
\begin{displaymath}
	I\ge
	\Cn
	\sum_{i=N+1}^{M-1}
		2^i\lambda_{2^i}
		\biggl(
			\sum_{j=i}^{M-1}
				a_{2^{j+1}}
				\sum_{\nu=2^j+1}^{2^{j+1}}\gamma_\nu
		\biggr)^p
	\ge
	\Cn
	\sum_{i=N+1}^{M-1}
		2^i\lambda_{2^i}
		\biggl(
			\sum_{j=i}^{M-1}
				a_{2^{j+1}} 2^j\gamma_{2^j}
		\biggr)^p.
\end{displaymath}
By applying Theorem~\ref{th:jensen} to this inequality
taking into account that $1\le p$,
then changing the order of summation,
we have
\begin{multline*}
	I\ge
	\lastC
	\sum_{i=N+1}^{M-1}
		2^i\lambda_{2^i}
		\sum_{j=i}^{M-1}
			a_{2^{j+1}}^p 2^{jp}\gamma_{2^j}^p
	\ge
	\lastC
	\sum_{j=N+1}^{M-1}
		a_{2^{j+1}}^p 2^{jp}\gamma_{2^j}^p
		\sum_{i=N+1}^j
			2^i\lambda_{2^i}\\
	\ge
	\lastC
	\sum_{j=N+1}^{M-1}
		a_{2^{j+1}}^p 2^{j(p+1)}\gamma_{2^j}^p\lambda_{2^j}.
\end{multline*}
Since $a_\nu\downarrow$,
taking into consideration that
$\seq{\gamma_\nu}{\nu=1}\infty$
and $\seq{\lambda_\mu}{\mu=1}\infty$
are quasi lacunary monotone
and quasi geometrically increasing sequences,
respectively,
we obtain
\begin{displaymath}
	I\ge
	\Cn
	\sum_{j=N+1}^{M-1}
		\sum_{\mu=2^{j+1}+1}^{2^{j+2}}
			a_\mu^p\gamma_\mu^p\mu^p\lambda_\mu.
\end{displaymath}
We rewrite the above inequality in the form
\begin{displaymath}
	I\ge
	\lastC
	\sum_{\mu=2^{N+2}+1}^{2^{M+1}}
		\lambda_\mu(\mu a_\mu\gamma_\mu)^p,
\end{displaymath}
and since $2^{N+2}<8n$, we obtain
\begin{displaymath}
	I\ge
	\lastC
	\sum_{\mu=8n}^n
		\lambda_\mu(\mu a_\mu\gamma_\mu)^p.
\end{displaymath}

Thus, we have proved inequality~\eqref{eq:lambda-p-ge-pge1}
assuming that $N+1\le M-1$.
In fact, for $m\ge16n$
we get $2^{N-1}<n\le\frac m{16}\le2^{M-3}$,
yielding that the condition $N+1\le M-1$ is satisfied.

In order to prove inequality~\eqref{eq:lambda-m-ge-pge1},
put
\begin{displaymath}
	J=
	\sum_{\mu=n}^m\lambda_\mu
		\biggl(
			\sum_{\nu=n}^\mu a_\nu\gamma_\nu
		\biggr)^p,
\end{displaymath}
in a similar manner,
but by making use of the fact that
$\seq{\lambda_\mu}{\mu=1}\infty$
is solely a quasi lacunary monotone sequence
(i.e.\ without a quasi geometrically increasing sequence assumption),
we obtain
\begin{displaymath}
	J\ge
	\Cn
	\sum_{\mu={2^{N+1}}}^{2^M-1}
		\lambda_\mu(\mu a_\mu\gamma_\mu)^p.
\end{displaymath}
Thus
\begin{equation}\label{eq:sum1-le-J}
	\sum_{\mu={2^{N+1}}}^{2^M-1}
		\lambda_\mu(\mu a_\mu\gamma_\mu)^p
	\le
	\Cn
	J.	
\end{equation}
Since $a_\nu\downarrow$
and $\seq{\lambda_\mu}{\mu=1}\infty$,
$\seq{\gamma_\nu}{\nu=1}\infty$
are lacunary monotone sequences,
we have
\begin{displaymath}
	J_1=
	\sum_{\mu={2^M}}^{2^{M+1}}
		\lambda_\mu(\mu a_\mu\gamma_\mu)^p
	\le
	a_{2^M}^p
	\sum_{\mu={2^M}}^{2^{M+1}}
		\lambda_\mu(\mu a_\mu\gamma_\mu)^p\\
	\le
	\Cn
	\sum_{\mu={2^{M-1}}}^{2^M-1}
		\lambda_\mu(\mu a_\mu\gamma_\mu)^p.
\end{displaymath}
Hence,
for $N+1\le M-1$ we obtain
\begin{displaymath}
	J_1\le
	\lastC
	\sum_{\mu={2^{N+1}}}^{2^M-1}
		\lambda_\mu(\mu a_\mu\gamma_\mu)^p.
\end{displaymath}
Thus, inequality~\eqref{eq:sum1-le-J} yields
\begin{displaymath}
	J_1\le\Cn J;
\end{displaymath}
or
\begin{equation}\label{eq:sum2-le-J}
	\sum_{\mu={2^M}}^{2^{M+1}}
		\lambda_\mu(\mu a_\mu\gamma_\mu)^p.
	\le
	\lastC
	J.
\end{equation}
Adding inequalities~\eqref{eq:sum1-le-J}
and~\eqref{eq:sum2-le-J} together,
we obtain
\begin{displaymath}
	\sum_{\mu={2^{N+1}}}^{2^{M+1}}
		\lambda_\mu(\mu a_\mu\gamma_\mu)^p
	\le
	\Cn
	\sum_{\mu=n}^m\lambda_\mu
		\biggl(
			\sum_{\nu=n}^\mu a_\nu\gamma_\nu
		\biggr)^p.
\end{displaymath}
Since $2^{N+1}\le4n$, $2^{M+1}>m$,
the above inequality implies the inequality~\eqref{eq:lambda-m-ge-pge1}.
This completes the proof of Theorem~\ref{th:lambda-ge-pge1}.

\section{Proof of Theorem~\ref{th:lambda-le-ple1}}

We prove the inequality~\eqref{eq:lambda-p-le-ple1}.
Let positive integers~$N$ and~$M$
be defined by the inequalities
$2^{N-1}<n\le2^N$
and $2^M\le m<2^{M+1}$.
This yields
\begin{multline*}
	I=
	\sum_{\mu=4n}^m\lambda_\mu
		\biggl(
			\sum_{\nu=\mu}^m a_\nu\gamma_\nu
		\biggr)^p
	\le
	\sum_{\mu=2^{N+1}+1}^{2^{M+1}}\lambda_\mu
		\biggl(
			\sum_{\nu=\mu}^{2^{M+1}} a_\nu\gamma_\nu
		\biggr)^p\\
	=
	\sum_{i=N+1}^M
		\sum_{\mu=2^i+1}^{2^{i+1}}
			\lambda_\mu
			\biggl(
				\sum_{\nu=\mu}^{2^{M+1}}
					a_\nu\gamma_\nu
			\biggr)^p.
\end{multline*}
By bounding the third sum from above,
taking into account that $\seq{\lambda_\mu}{\mu=1}\infty$
is a quasi lacunary monotone sequence,
we get
\begin{displaymath}
	I\le
	\sum_{i=N+1}^M
		\biggl(
			\sum_{\nu=2^i+1}^{2^{M+1}}
				a_\nu\gamma_\nu
		\biggr)^p
		\sum_{\mu=2^i+1}^{2^{i+1}}
			\lambda_\mu
	\le
	\Cn
	\sum_{i=N+1}^M
		2^i\lambda_{2^i}
		\biggl(
			\sum_{\nu=2^i+1}^{2^{M+1}}
				a_\nu\gamma_\nu
		\biggr)^p,
\end{displaymath}
where positive constant~$\lastC$
depends only on the sequence $\seq{\lambda_\mu}{\mu=1}\infty$.
Now, we split the second sum into blocks of length $2^j$
and taking into consideration the fact that $a_\nu\downarrow$,
and $\seq{\gamma_\nu}{\nu=1}\infty$
is a lacunary monotone sequences,
we have
\begin{displaymath}
	I\le
	\lastC
	\sum_{i=N+1}^M
		2^i\lambda_{2^i}
		\biggl(
			\sum_{j=i}^M
				a_{2^j}
				\sum_{\nu=2^j+1}^{2^{j+1}}
					\gamma_\nu
		\biggr)^p
	\le
	\Cn
	\sum_{i=N+1}^M
		2^i\lambda_{2^i}
		\biggl(
			\sum_{j=i}^M
				a_{2^j} 2^j\gamma_{2^j}
		\biggr)^p.
\end{displaymath}

By applying Theorem~\ref{th:jensen}
and then changing the order of summation,
we obtain
\begin{displaymath}
	I\le
	\lastC
	\sum_{i=N+1}^M
		2^i\lambda_{2^i}
		\sum_{j=i}^M
			a_{2^j}^p 2^{jp}\gamma_{2^j}^p
	\le
	\lastC
	\sum_{j=N+1}^M
		a_{2^j}^p 2^{jp}\gamma_{2^j}^p
		\sum_{i=N+1}^j
			2^i\lambda_{2^i}.
\end{displaymath}
Further,
the fact that
$\seq{\lambda_\mu}{\mu=1}\infty$
is a quasi geometrically increasing sequence
yields
\begin{displaymath}
	I\le
	\Cn
	\sum_{j=N+1}^M
		a_{2^j}^p 2^{jp}\gamma_{2^j}^p
		2^j\lambda_{2^j}.
\end{displaymath}
Since $a_\nu\downarrow$,
taking into consideration that
$\seq{\gamma_\nu}{\nu=1}\infty$
and $\seq{\lambda_\mu}{\mu=1}\infty$
are quasi lacunary monotone sequences,
we get
\begin{displaymath}
	a_{2^j}^p 2^{j(p+1)}\gamma_{2^j}^p\lambda_{2^j}
	\le
	\Cn
	\sum_{\mu=2^{j-1}+1}^{2^j}
		a_\mu^p\gamma_\mu^p\mu^p\lambda_\mu.
	\end{displaymath}
Therefore,
\begin{displaymath}
	I\le
	\Cn
	\sum_{j=N+1}^M
		\sum_{\mu=2^{j-1}+1}^{2^j}
			\lambda_\mu(\mu a_\mu\gamma_\mu)^p
	=
	\lastC
	\sum_{\mu=2^N+1}^{2^M}
		\lambda_\mu(\mu a_\mu\gamma_\mu)^p.
\end{displaymath}
where positive constant~$\lastC$
does not depend on~$N$ and~$M$.
Since $2^M\le m$ and $2^N+1>m$,
we obtain
\begin{displaymath}
	I\le
	\lastC
	\sum_{\mu=n}^m
		\lambda_\mu(\mu a_\mu\gamma_\mu)^p,
\end{displaymath}
which proves inequality~\eqref{eq:lambda-p-le-ple1}.

Inequality~\eqref{eq:lambda-m-le-ple1}
can be proved in an analogous way,
but without a quasi geometrically increasing assumption
for the sequence
$\seq{2^\mu\lambda_{2^\mu}}{\mu=1}\infty$.

\section{Inequalities for non-negative decreasing sequences}

For $\alpha, \lambda\in\numR$
put
\begin{gather*}
	\lambda_\mu=\mu^{\alpha - 1}
		\quad (\mu=1,2,\dots),\\
	\gamma_\nu=\nu^\lambda
		\quad (\nu=1,2,\dots).
\end{gather*}
Obviously,
the obtained sequences $\seq{\mu^{\alpha - 1}}{\mu=1}\infty$
and $\seq{\nu^\lambda}{\nu=1}\infty$
are both quasi lacunary monotone sequences.

If, in addition, $\alpha>0$,
then $\seq{2^\mu 2^{\mu(\alpha - 1)}}{\mu=1}\infty$
is a quasi geometrically increasing sequence.

By applying Theorems~\ref{th:lambda-ge-pge1}
and~\ref{th:lambda-le-ple1}
for such sequences $\seq{\lambda_\mu}{\mu=1}\infty$
and $\seq{\gamma_\nu}{\nu=1}\infty$,
we deduce the following
$l_p$-type inequalities for non-negative decreasing sequences,
which are converse to inequalities~\eqref{eq:p-le-pge1},
\eqref{eq:m-le-pge1}, \eqref{eq:p-ge-ple1}
and~\eqref{eq:m-ge-ple1}.

\begin{theorem}\label{th:ge-pge1}
	Let a sequence $\seq{a_\nu}{\nu=1}\infty$ be such that $a_\nu\downarrow$,
	and $\alpha>0$, $\lambda\in\numR$,
	$m$ and~$n$ positive integers such that $m\ge16n$.
	If $p\ge1$,
	then the following inequalities hold
	\begin{gather*}
		\sum_{\mu=n}^m \mu^{\alpha-1}
			\biggl(
				\sum_{\nu=\mu}^m a_\nu \nu^\lambda
			\biggr)^p
		\ge C
		\sum_{\mu=8n}^m \mu^{\alpha-1}
			(a_\mu \mu^{\lambda+1})^p,
		\\
		\sum_{\mu=n}^m \mu^{-\alpha-1}
			\biggl(
				\sum_{\nu=n}^\mu a_\nu \nu^\lambda
			\biggr)^p
		\ge C
		\sum_{\mu=4n}^m \mu^{-\alpha-1}
			(a_\mu \mu^{\lambda+1})^p.
	\end{gather*}
	Hereafter
	$C$~denotes positive constant
	depending only on~$\alpha$, $\lambda$ and~$p$,
	and not depending on~$m$, $n$
	and the sequence $\seq{a_\nu}{\nu=1}\infty$.
\end{theorem}
                           
\begin{theorem}\label{th:le-ple1}
	Let a sequence $\seq{a_\nu}{\nu=1}\infty$ be such that $a_\nu\downarrow$,
	and $\alpha>0$, $\lambda\in\numR$,
	$m$ and~$n$ positive integers such that $m\ge4n$.
	If $0<p\le1$,
	then the following inequalities hold
	\begin{gather*}
		\sum_{\mu=4n}^m \mu^{\alpha-1}
			\biggl(
				\sum_{\nu=\mu}^m a_\nu \nu^\lambda
			\biggr)^p
		\le C
		\sum_{\mu=n}^m \mu^{\alpha-1}
			(a_\mu \mu^{\lambda+1})^p,
		\\
		\sum_{\mu=4n}^m \mu^{-\alpha-1}
			\biggl(
				\sum_{\nu=4n}^\mu a_\nu \nu^\lambda
			\biggr)^p
		\le C
		\sum_{\mu=n}^m \mu^{-\alpha-1}
			(a_\mu \mu^{\lambda+1})^p.
	\end{gather*}
\end{theorem}

Note that Theorems~\ref{th:ge-pge1} and~\ref{th:le-ple1}
given above imply several inequalities
proved earlier~\cite{hardy-l-p:inequalities,
	konyushkov:mat-58, vukolova:vestnik-84%
}.

Namely the following Corollaries are simple consequences
of these theorems and the inequalities~\eqref{eq:p-le-pge1},
\eqref{eq:m-le-pge1}, \eqref{eq:p-ge-ple1} 
and~\eqref{eq:m-ge-ple1}.

\begin{corollary}\label{cr:ge}
	Let a sequence $\seq{a_\nu}{\nu=1}\infty$ be such that $a_\nu\downarrow$,
	$\alpha>0$, $\lambda\in\numR$,
	and~$n$ a positive integer.
	If $p>0$,
	then the following inequalities hold
	\begin{gather*}
		\sum_{\mu=1}^n \mu^{\alpha-1}
			\biggl(
				\sum_{\nu=\mu}^n a_\nu \nu^\lambda
			\biggr)^p
		\ge C
		\sum_{\mu=1}^n \mu^{\alpha-1}
			(a_\mu \mu^{\lambda+1})^p,
		\\
		\sum_{\mu=1}^n \mu^{-\alpha-1}
			\biggl(
				\sum_{\nu=1}^\mu a_\nu \nu^\lambda
			\biggr)^p
		\ge C
		\sum_{\mu=1}^n \mu^{-\alpha-1}
			(a_\mu \mu^{\lambda+1})^p.
	\end{gather*}
\end{corollary}

\begin{corollary}\label{cr:asymp}
	Let a sequence $\seq{a_\nu}{\nu=1}\infty$ be such that $a_\nu\downarrow$,
	$\alpha>0$, $\lambda\in\numR$,
	and~$n$ a positive integer.
	If $p\ge1$,
	then the following asymptotic equivalences hold
	\begin{gather*}
		\sum_{\mu=1}^n \mu^{\alpha-1}
			\biggl(
				\sum_{\nu=\mu}^n a_\nu \nu^\lambda
			\biggr)^p
		\asymp
		\sum_{\mu=1}^n \mu^{\alpha-1}
			(a_\mu \mu^{\lambda+1})^p,
		\\
		\sum_{\mu=1}^n \mu^{-\alpha-1}
			\biggl(
				\sum_{\nu=1}^\mu a_\nu \nu^\lambda
			\biggr)^p
		\asymp
		\sum_{\mu=1}^n \mu^{-\alpha-1}
			(a_\mu \mu^{\lambda+1})^p.
	\end{gather*}
\end{corollary}

\section*{Acknowledgement}
The authors would like to thank the referee of the paper
for providing insightful comments
which have contributed substantially
to the overall quality of the article.

\bibliographystyle{hplain}
\bibliography{maths}

\end{document}